\documentclass[12pt,reqno]{amsart}

\usepackage[margin=1in]{geometry} 

\usepackage{times}
\usepackage[english]{babel}
\usepackage[utf8]{inputenc}
\usepackage{csquotes}
\usepackage{cancel}
\usepackage{hyperref}
\usepackage{tikz}
\usetikzlibrary{patterns,calc}

\usepackage{amsmath,amssymb,amsthm,bbm}
\usepackage{cases}
\usepackage{amsfonts}
\usepackage{mathtools}
\usepackage{esint}
\usepackage{tabulary}
\usepackage{bm}
\usepackage{graphicx, color}
\usepackage{svg}
\usepackage{enumitem}
\usepackage{cancel}
\usepackage{parskip}
\newcommand{\dx}{\,{\rm d} x}

\newcommand{\dt}{\,{\rm d} t}
\newcommand{\ds}{\,{\rm d} s}

\newcommand{\numberset}{\mathbb}
\newcommand{\R}{\numberset{R}}

\newcommand{\N}{\numberset{N}}
\newcommand{\C}{\numberset{C}}

\newcommand{\Rz}{\R}
\newcommand{\SO}{{\rm SO}}

\newcommand{\GL}{{\rm GL}}
\newcommand{\J}{V^{J}}

\newcommand{\RR}{\mathcal{R}}
\newcommand{\HH}{\mathcal{H}}
\newcommand{\W}{\mathcal{W}}
\newcommand{\A}{\mathcal{A}}

\newcommand{\DD}{\mathbb{D}}

\DeclareMathOperator*{\argmin}{arg\,min}

\newcommand{\weak}{\rightharpoonup}

\newtheorem{theorem}{Theorem}[section]
\newtheorem{definition}{Definition}[section]
\newtheorem{proposition}{Proposition}[section]

\theoremstyle{definition}

\newcommand{\EEE}{\color{black}}

\title[ Existence for accreting viscoelastic solids 
]{ Existence for accreting viscoelastic solids at large
  strains
}

\author[A. Chiesa] {Andrea Chiesa} 
\address[Andrea Chiesa]{University of Vienna, Faculty of Mathematics  and Vienna School of Mathematics,
                Oskar-Morgenstern-Platz 1, A-1090 Vienna, Austria}
\email{andrea.chiesa@univie.ac.at}
	\urladdr{http://www.mat.univie.ac.at/$\sim$achiesa}
\author[U. Stefanelli]{Ulisse Stefanelli} 
	\address[Ulisse Stefanelli]{University of
		Vienna, Faculty of Mathematics,
                Oskar-Morgenstern-Platz 1, A-1090 Vienna, Austria, 
		University of Vienna, and Vienna Research Platform on Accelerating
		Photoreaction Discovery, W\"ahringerstra\ss e 17, 1090
                Wien, Austria}
	\email{ulisse.stefanelli@univie.ac.at}
	\urladdr{http://www.mat.univie.ac.at/$\sim$stefanelli}


\allowdisplaybreaks
\begin{document}

\begin{abstract}
   By revisiting  a model proposed in
\cite{TruskinovskyZurlo2}, we  address   the accretive growth of a viscoelastic 
solid at large strains. The
accreted material is assumed to  accumulate at the boundary of the body  in an
unstressed state.  The growth process is driven by the deformation
state of the solid.   The \EEE  progressive build-up
of incompatible strains in the material   is modeled \EEE by considering an additional
backstrain. The model is regularized by postulating  the presence
of  a fictitious compliant material  surrounding the accreting body. 
We   show the existence of solutions to the coupled accretion and viscoelastic equilibrium problem.  
\end{abstract}

\maketitle

\section{Introduction}

Growth is a fundamental process in all biological systems, as
well as in a
variety of natural, technological, and social ones. Among the many different
dynamics, {\it accretive growth} occurs when growth is realized via a
progressive accumulation, addition, or layering   of material \EEE
{\it at the boundary} of the system. This paradigm
is of paramount relevance in numerous situations. The formation of horns, teeth, and seashells
\cite{Moulton,calcific1, calcific2}, coral reefs \cite{Kaandorp}, bacterial
colonies \cite{Hoppensteadt},  trees
\cite{plants,trees}, and cell motility due to actin growth
\cite{hodge} are biological
examples of accretive growth. In geophysics, sedimentation and glacier
formations   are also accretive processes, as is planet
formation \cite{planet}.
 Furthermore,  accretive growth is a key aspect in many
technological applications,  including metal solidification
\cite{solidification1}, crystal growth \cite{crystallization0,crystallization2},  
additive manufacturing~\cite{3Dprint1,horn,3Dprint2}, layering, coating, and masonry,
just to mention a few.

In this note, we consider the evolution in time of a viscoelastic solid under accretive
growth. Correspondingly, the reference configuration of the body
$\Omega(t)\subset \Rz^d\ (d \geq 2)$ is time-dependent and the deformation
$y(t,\cdot):\Omega(t) \to \Rz^d$ of the solid is defined on a time-dependent
domain. 
Although in some cases the map
$t \in [0,T] \mapsto \Omega(t)$ can be rightfully assumed to be given
(this is for instance the case for 3D printing) the evolution in time of the reference configuration
is not a-priori known in general, but is rather influenced by the mechanical
process. In this paper, we assume $t \mapsto \Omega(t)$ to be unknown
and we tackle its specification by  adopting a {\it level-set approach}
  \cite{Evans-Spruck,Souganidis} \EEE and setting
 \begin{equation}
   \Omega(t):=\{x\in \Rz^d \:|\: \theta(x)<t\}.\label{eq:zero}
 \end{equation}
The map $x \mapsto \theta(x)\in [0,\infty)$ is called
\emph{time-of-attachment} function: the value  $\theta(x)$ corresponds
to  the instant in time at which the point $x \in
\Rz^d$ the accreting body   reaches the point $x \in \Rz^d$. \EEE We assume that accretion 
 occurs at a positive {\it growth rate} $\gamma(\cdot)>0$ and in the outward pointing normal
 direction to the body.

 This implies that the
time-of-attachment function $\theta$ solves the following
external problem for  
the {\it generalized eikonal equation}  
\cite{Soner,DavoliNikStefanelliTomassetti}
\begin{align}\label{eq growth intro paper 2}
  &\gamma(\nabla y(\theta(x),x))|\nabla \theta(x)|=1 \quad \text{for}
    \ x\in \Omega(T)\setminus \overline{\Omega_0},\\
  &\theta(x)=0  \quad \text{for} \  x\in  \Omega_0\label{eq growth intro paper 22}.
\end{align} 
Here, $\Omega_0\subset \R^d$ is the given initial reference
configuration of the accreting   solid. \EEE As growth is often driven by the mechanical state of the
accreting body \cite{Goriely}, we let the
growth rate $\gamma$ by setting 
 $\gamma=\gamma(\nabla y(\theta(x),x))$. Specifically, $\gamma\colon
\R^{d\times d}\rightarrow [0,\infty)$ is assumed to be Lipschitz
continuous and
such that $c_\gamma\leq \gamma(\cdot) \leq C_\gamma$ for some
$0<c_\gamma\leq C_\gamma$.

In order to track the progressive accumulation of
incompatibilities due to growth \cite{solidification1,
  calcific1,TruskinovskyZurlo2}, we follow the classical approach of
finite plasticity \cite{Kroener,Lee} and postulate the {\it multiplicative
  decomposition} 
$$\nabla y = F_{\rm e} A$$
where $F_{\rm e}\in \GL_+(d)$ corresponds to the elastic part of the
deformation gradient, whereas  $A\in \GL_+(d)$ is the backstrain originated by the
build-up of incompatibilities \cite{Goriely}   during growth. The \EEE elastic energy
density $W\colon\Rz^{d\times d}\to [  0,\infty)$ of the accreting medium is assumed to be depending on $F_{\rm e} =
\nabla y A^{-1}$.  In order to specify a constitutive relation for
$A$, we  follow \cite{TruskinovskyZurlo2},  see also 
\cite{DavoliNikStefanelliTomassetti,TruskinovskyZurlo1}, and assume that  the material
is added to the accreting body in an unstressed state. Specifically, we
assume $W$ to be minimized at the identity matrix $I$ and $F_{\rm
  e}(t,x)=I$ at the accreting front, i.e., for $(t,x)=(\theta(x),x)$. This
entails the constitutive relation 
\begin{equation}
   A(x)  = \nabla y(\theta(x),x) \quad \text{for} \ x\in 
  \Omega(T)\setminus \overline {\Omega_0}. \label{eq:A}
\end{equation} 
Note that the backstrain $A$ is independent of time. In particular, as
the boundary 
$\partial \Omega(t)$ reaches point $x$, the value of $A(x)$ is stored
according to \eqref{eq:A}. This reflects the intuition that growth-driven
incompatibilities are recorded in the body along the process.

The viscoelastic evolution of the accreting solid is determined by the
equilibrium system
\begin{align}\label{eq:y growth 0}
&{-}\operatorname{div}\!\left( {\rm D}W ( \nabla
    yA^{-1})A^{-\top} + {\rm D}\J( \nabla
    y) + \partial_{\nabla \dot y}R(\nabla y{,}\nabla \dot
                                     y) - \operatorname{div}{\rm D}
                                     H(\nabla^2y)\right)   = f
\end{align}
to be solved in the noncylindrical domain $\cup_{t\in [0,T]}\{t\}
\times \Omega(t)$.
 Here,   $\J\colon \R^{d\times d}\rightarrow [0,\infty]$ is an
additional elastic energy term 
 specifically penalizing  self-interpenetration of matter, i.e.,
$\J(F)\rightarrow \infty$ as $\det F\rightarrow 0^+$ and
$\J(F)<\infty$ if and only if $\det F>0$.  Moreover,
$R:\GL_+(d)\times \Rz^{d\times d}   \rightarrow [0,\infty  ) \EEE$ is the
instantaneous dissipation potential, which is assumed to be quadratic
in $\nabla \dot y ^\top\nabla y +\nabla y ^\top\nabla \dot y $ where the dot
denotes the partial time derivative.  The term $H\colon
\R^{d\times d\times d}\rightarrow [0,\infty)$ qualifies the accreting solid
as a {\it second-grade,  nonsimple} material  
and $f=f(t,x)$ is  and external-force density.

\begin{figure}[h]
  \centering
  \pgfdeclareimage[width=150mm]{figure}{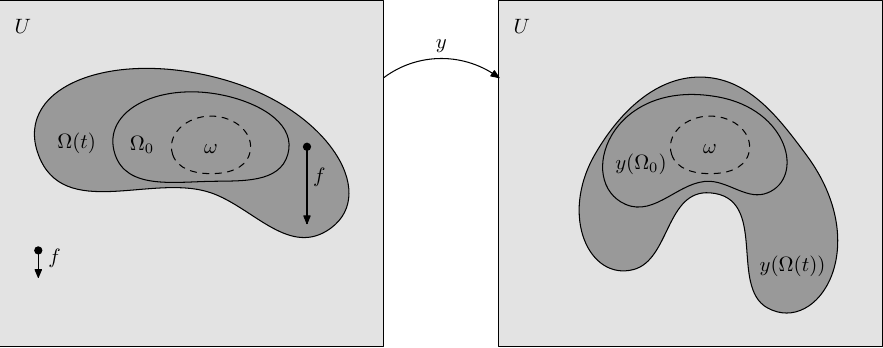}
  \pgfuseimage{figure}
\caption{Illustration of the notation in the reference setting  (left) and in the deformed one
(right), with the boundary and docking conditions \eqref{eq:cond2}--\eqref{attachment accretive 2 intro}.}
  \label{Fig:1}
\end{figure}

The model \eqref{eq:zero}--\eqref{eq:y growth 0} is of
{\it free-boundary} type, as equations are posed on the unknown sets
$\Omega(t)$ from \eqref{eq:zero}. 
  This creates significant difficulties for the analysis, forcing us to reduce the model to a fixed, ambient setting, see Figure \ref{Fig:1}. In
particular, we ask that $\Omega(t) \subset U $ for all
$t\in [0,T]$, for a fixed, open, and bounded
container $U \subset \Rz^d$.
Problem \eqref{eq growth
  intro paper 2}--\eqref{eq growth intro paper 22} can be reduced to
the fixed-boundary setting by considering
\begin{align}\label{eq growth intro paper 20}
  &\gamma(\nabla y(\theta(x)\wedge T,x))|\nabla \theta(x)|=1 \quad \text{for}
    \ x\in U\setminus \overline{\Omega_0}\\
  &\theta(x)=0  \quad \text{for} \  x\in  \Omega_0\label{eq growth intro paper 220},
\end{align}
instead of \eqref{eq growth intro paper 2}--\eqref{eq growth intro
  paper 20}. This modification is actually immaterial as   we will
check that \EEE 
the restriction to 
$\Omega(T)=\{x\in U \ | \ \theta(x)<T\}$ of a  solution $\theta$ to
\eqref{eq growth intro paper 22}--\eqref{eq growth intro paper 220}
solves \eqref{eq growth intro paper 2}--\eqref{eq growth intro paper
  22}, as well.  Note nonetheless that \eqref{eq growth intro
  paper 20} requires to introduce the minimum $\theta(x)\wedge
T$, as $\theta(x)>T$ at some points of $U$. 

\vspace{5mm}

\vspace{5mm}

In order to reduce the mechanical problem \eqref{eq:A}--\eqref{eq:y
  growth 0} to the fixed container, we assume that the complement
$U\setminus \overline{\Omega(t)}$ to the accreting body is filled by a second,
{\it fictitious} medium of viscoelastic type. By letting
such fictitious medium be very compliant, one expects that the
behavior of the regularized model to approximate the original,
free-boundary one.
Correspondingly, the constitutive equation
\eqref{eq:A} may be extended to $U$ by posing
\begin{equation}
  \label{eq:A0}
  A(  t, \EEE x)=
  \left\{
    \begin{array}{ll}
      A_0(x)  \quad&\text{for} \ x \in   \overline{\Omega_0},\\
    \nabla y(\theta(x),x) \quad&\text{for} \ x \in \Omega(  t \EEE)\setminus \overline{\Omega_0},\\
    I \quad&\text{for} \ x \in {U}\setminus \overline{\Omega(  t \EEE)},
    \end{array}
    \right.
  \end{equation}
  where we have additionally specified the initial backstrain $A_0$ on $\Omega_0$
  and defined $A=I$ in the region $ {U}\setminus \overline{\Omega( 
    t \EEE)}$
  which is not accessed by the accreting medium   at time $t$, \EEE
  by simplicity.   Compared with \eqref{eq:A}, the backstrain
  $A(t,x)$ is now depending on $t$, as well. This is an artifact of
  the extension of the model to the container $U$. \EEE

Eventually, the viscoelastic equilibrium \eqref{eq:y growth 0} also
need to be extended to the whole container. 
At all $(t,x)\in [0,T]\times U$, we may distinguish the
accreting medium and the fictitious one by the sign of
$\theta(x)-t$. More precisely, if $\theta(x)-t\leq0$ we have that $x\in
\overline{\Omega(t)}$ (accreting medium), whereas $\theta(x)-t> 0$ implies
that $x\in U \setminus \overline{\Omega(t)}$ (fictitious
material). The sharp transition in the material parameters between the
two media will be described by a function  $(x,t)\mapsto
h(\theta(x)-t)$ with $h(\theta(x)-t)=1$
for $\theta(x)-t \leq 0$ and  $h(\theta(x)-t)=  \delta\EEE$ 
for $\theta(x)-t > 0$ with
$\delta  \in (0,1) \EEE $ very small. By using $h$, the elastic energy density and
the instantaneous
dissipation density of combined medium are considered to be $h(\theta(x)-t)W(\nabla y A^{-1})$
and $h(\theta(x)-t)R(\nabla y, \nabla \dot y)$, respectively, 
and the  viscoelastic equilibrium  system 
  takes the form \EEE   
\begin{align}
&-\operatorname{div}\!\left(h(\theta(x){-}t){\rm D}W ( \nabla
    y(t,x)A^{-1}(  t, \EEE x))A^{-\top}(  t, \EEE x) +  {\rm D}  \J( \nabla
                                     y(t,x))\right)\notag\\
  &\quad -\operatorname{div}\!\left(h(\theta(x){-}t)\partial_{\nabla \dot y}R(\nabla y(t,x),\nabla \dot y(t,x)){-}\!\operatorname{div}{\rm D}  H(\nabla^2y(t,x))\right)\notag\\
  &\quad  =h(\theta(x){-}t)f(t,x) \quad \text{for} \ (t,x)\in
    (0,T)\times U. \label{eq:y growth 2}
\end{align} 
Notice that the second-order-potential density $H$ and the term $\J$
are assumed to be the same for both the  accreting and the fictitious
media. On the other hand, the force density $h(\theta(x){-}t)f(t{,}x)$
  distinguishes \EEE between the accreting and the fictitious
medium, as it would be the case for gravity for different densities.  


 The aim of this paper is to show that model \eqref{eq growth
  intro paper 20}--\eqref{eq:y growth 2} admits 
solutions. To this end,
system \eqref{eq growth intro paper 20}--\eqref{eq:y growth 2} is
complemented by the boundary and
initial conditions 
\begin{align}
  &{\rm D}H(\nabla^2y) {:}  (\nu\otimes \nu)=0 \ \  \text{on}  \ \
      (0,T) \EEE \times\partial U,\label{eq:cond1}\\
   &y={\rm id}\ \ \text{on} \ \    (0,T) \EEE\times\partial U, \label{eq:cond2}\\
  \label{attachment accretive 2 intro}
  &y = {\rm id} \ \ \text{on}  \ \   (0,T) \EEE\times \omega,\\
  & y(0,\cdot )= y_0 \ \ \text{on}  \ \ U. \label{eq:cond3}
\end{align}
 
Specifically, at the  boundary $\partial U$ of the
container $U$ we prescribe the homogeneous Neumann and Dirichlet conditions
\eqref{eq:cond1}--\eqref{eq:cond2}. These choices are motivated by
simplicity,   reflecting \EEE the mere instrumental role of the container $U$
in the model. In fact, other options would be viable, as well.   
On the other hand,  the \emph{anchoring condition} \eqref{attachment
  accretive 2 intro} fixes  the position of the  body at
 a portion $\omega\subset \subset \Omega_0$ of the
starting configuration $\Omega_0$  see 
\cite{DavoliNikStefanelliTomassetti}.  
Independently of \eqref{eq:cond1}--\eqref{eq:cond2}, this anchoring
condition will allow to use a Poincar\'e-type inequality, cf. \eqref{Poincare}, which turns out to
be crucial for the analysis. 

 Our main result consists in proving existence of solutions $(\theta,y)$
to the fully coupled system  %
\eqref{eq growth intro paper 20}--\eqref{eq:cond3}, see Definition
\ref{def:weak sol A} and  Theorem \ref{Thm:existence sol A}. 
More precisely, we   find \EEE $\theta \in C(\overline{U})$  
solving \EEE  the
external problem for the generalized eikonal equation \eqref{eq growth
  intro paper 20}--\eqref{eq growth intro paper 220} in the viscosity
sense \cite{Viscosity sol user's guide},   and \EEE 
$y \in L^\infty(0,T;W^{2,p}(U;\Rz^d))\cap H^1(0,T;H^1(U;\Rz^d) )$
  satisfying \EEE  the viscoelastic equilibrium system \eqref{eq:y growth
  2}--\eqref{eq:cond3} in
 the weak sense,   with \EEE the constitutive relation \eqref{eq:A0} 
 pointwise fulfilled. 

  Before proceeding, let us mention that the applied literature on
the mechanics of growth is well developed. Comprehensive accounts can
be found in the monographs \cite{Goriely,calcific2}, as well as in
\cite{Ganghoffer2,finitestrainlit,WalkingtonDayal1,WalkingtonDayal2,Sozio}, among
many others. In comparison, rigorous mathematical results mechanical
growth models are just a few.  For instance,  an elastic
growing body and the coupled dynamics of the morphogen is studied in
\cite{existence1}. In the setting of {\it
  bulk} growth, where growth is realized by addition of material in
the bulk of the solid, existence results are available in both 
one 
\cite{Bangert} and three
\cite{existence2,Ganghoffer} space dimensions.

In
the specific setting of accretive growth, which is the focus of this paper, a first existence result in
the context of linearized elasticity has been obtained in
\cite{DavoliNikStefanelliTomassetti}. There, the constitutive relation
\eqref{eq growth intro paper 2} cannot be directly considered  due to the
limited regularity of the solution \cite{john}, and  an
additional regularization via a mollification is introduced, which can be interpreted as a
diffused-interface, phase-field approximation. By neglecting the
backstress $A$, the viscoelastic growth model  \eqref{eq growth intro paper 2}--\eqref{eq:y growth 2} 
has been considered in \cite{viscoelastic_growth}, where an existence
result is obtained, both in the mollified setting (as above) and the
limiting sharp-interface
case.

\EEE

  The paper is structured as follows. \EEE
In Section \ref{sec: notation}, we  specify notation and the
assumptions on the ingredients of the   model.  The definition
of   solution to the problem \EEE is detailed in Section \ref{sec: main
  results A}, where the main existence result, Theorem
\ref{Thm:existence sol A}, is also stated. Section \ref{sec:proof A}
is then devoted  to the proof of Theorem \ref{Thm:existence sol
  A}. The  proof  strategy  is iterative. At first, we
show that, for given $y$, we can find a viscosity solution $\theta$ to
\eqref{eq growth intro paper 20}--\eqref{eq growth intro paper 220}, see
Proposition \ref{prop:mennucci A}. Then, in Proposition \ref{prop:
  existence y a} we check the existence of $y$ solving
\eqref{eq:A0}--\eqref{eq:cond3} for given $\theta$. Eventually, we combine these
results   and \EEE  iteratively  define   a sequence
$(y_k,\theta_k)_{k\in \N}$   converging, \EEE up
to subsequences, to a solution to the fully coupled problem.

\section{Notation and setting}\label{sec: notation}

 We devote this section to specifying notation and   assumptions.

\subsection{Notation}
In what follows, we denote by $\R^{d\times d}$ the Euclidean space of
$d{\times} d$ real matrices,  $d\ge 2$,  by  $\R^{d\times d}_{\rm sym}$
the  subspace   of symmetric matrices, and by $I$ the identity
matrix. Given $A\in \R^{d\times d}$, we
indicate by $A^\top$ its  transpose and by $
    |A|^2 := A{:}A$ its Frobenius
  norm,
  where the contraction product between  two  matrices  $A,\,B\in
  \R^{d\times d}$  is defined as $A {:}
   B  :=  A_{ij}B_{ij}$  (we use the summation convention over
   repeated indices).  Analogously, let $\R^{d\times d\times d}$ be the set of real $3$-tensors, and define their contraction product as $A \smash{\vdots}
   B  := A_{ijk}B_{ijk}$ for $A,B\in \R^{d\times d\times d}$.
  Moreover, given a real 4-tensor
  $\C\in\R^{d\times d\times d\times d}$ and the matrix $A\in
  \R^{d\times d}$ we indicate by $\C{:}A\in \R^{d\times d}$ and $ A{:}\C\in \R^{d\times d}$ the matrices
  given in components by $(\C{:}A)_{ij} = \C_{ijk\ell}A_{k\ell}$ and
  $(A{:}\C)_{ij} = A_{k\ell}\C_{k\ell i j}$, respectively.   
  We shall use the following matrix sets
  \begin{align*}
    \SO(d)&\coloneqq \{ A \in \R^{d\times d} \;|\;\det  A=1,\,  A A^{\top}
   =I\},\\
   \GL_{+}(d)&\coloneqq\{ A  \in \R^{d\times d} \;|\; \det  A >0 \}.
  \nonumber
  \end{align*} 
   The scalar product of two vectors $a,\, b\in \R^d$ is
  classically indicated by $a{\cdot}b$. 
  The symbol $B_R\subset \R^{d}$ denotes the open ball of radius $R>0$
  and center $0\in\R^{d}$, $|E|$ indicates the Lebesgue measure
  of the Lebesgue-measurable set $E\subset \Rz^d$,  and
  $\mathbbm{1}_E$ is the corresponding characteristic function,
  namely, $\mathbbm{1}_E(x)=1$ for $x\in E$ and $\mathbbm{1}_E(x )=0$
  otherwise. For $E \subset \Rz^d$ nonempty and $x\in \Rz^d$ we define
  ${\rm dist}(x,E):=\inf_{e\in E}|x{-}e|$. We denote by
  $\mathcal{H}^{d-1}$ the $(d{-}1)$-dimensional Hausdorff measure 
  and by $\mathcal{L}^{d+1}$ the Lebesgue measure in $\Rz^{d+1}$. 
  
   In the following, we indicate by  $c$ a generic positive constant, possibly depending
on data but independent of  the time discretization step
$\tau$,  to be used in the proof of Proposition \ref{prop:
  existence y a}.  Specifically, $c$ may depend on $\delta>0$,  defined in \eqref{eq
  h} below. Note  that the value of $c$ may change from line to line.

\subsection{Setting}
 We start by posing the following. 
\begin{enumerate}[label=(H\arabic*)]
  \item\label{W0 A} 
Let $T>0$ be a fixed final time, $U \subset \R^d$ be nonempty,
open, connected, bounded, and Lipschitz, $  \Omega_0\subset \subset
U$   and $  \omega \subset \subset \Omega_0 $   be 
nonempty and open, and $p >d$. 
\end{enumerate}
We define $Q\coloneqq (0,T)\times U$.
 
\subsubsection{Admissible deformations}
The set of admissible deformations is  given  as 
\begin{align*}
  \A \coloneqq \Big\{y\in W^{2,p}_{\omega}(U;\R^d) \mid 
  \nabla y\in \GL_+(d) \text{ a.e. in } 
     U \Big\},
\end{align*}
where 
\begin{equation*}
  W^{2,p}_{\omega}(U;\R^d)\coloneqq \{y\in W^{2,p}(U;\R^d)\mid y\equiv {\rm id} \text{ on } \omega \cup \partial U\}.
\end{equation*}
Deformations $y$ are locally invertible and orientation
preserving. Moreover, they  satisfy the  anchoring condition
$y\equiv {\rm id}$ in $\omega$ for almost every $t\in (0,T)$. In
particular,  this   condition entails the validity of the following Poincar\'e-type inequality 
\begin{equation}\label{Poincare}
  \|y\|_{W^{2,p}(U;\R^d)}\leq c\left(1+\|\nabla^2 y\|_{L^p(U;\R^{d\times d\times d})}\right)\quad \forall y\in W^{2,p}_\omega(U;\R^d).
 \end{equation}

\subsubsection{ Energy}

The elastic energy density $W\colon \R^{d\times d}\rightarrow
[0,\infty)$ of the  accreting medium is asked to satisfy  \begin{enumerate}[label=(H\arabic*)]
  \setcounter{enumi}{1}
  \item\label{W1 A} $W\in C^1(\R^{d\times d})$;
  \item\label{W2 A} there exists $c_W>0$ such that 
  \begin{equation*}
    0=W(I)\leq   W( F) \leq \frac{1}{c_W} (|F|^{p}+1) \quad  
    \forall \EEE F\in \R^{d\times d}.
  \end{equation*}
\end{enumerate} 
 Although it is  not strictly needed for the analysis, we additionally assume
\begin{enumerate}[label=(H\arabic*)]
  \setcounter{enumi}{2}
  \item\label{W3 A} frame indifference: \   $W (Q F)=W (F)$ for all $F \in \R^{d\times d}$ and $Q \in \SO(d)$;
  \item\label{W4 A} isotropy: \   $W ( F Q)=W (F)$ for all $F \in \R^{d\times d}$ and $Q \in \SO(d)$.
\end{enumerate}
We remark that both \ref{W3 A} and \ref{W4 A} are required for the
model to be frame indifferent.  Indeed, taking \eqref{eq:A0} into
account, for all rotation $Q\in \SO(d)$ we have 
\begin{equation*} W(Q\nabla y(t,x)\nabla y^{-1}(\theta(x),x)Q^\top)=W(\nabla y(t,x)\nabla y^{-1}(\theta(x),x)).
\end{equation*} 

The density $\J\colon\GL_+(d)\rightarrow[0,\infty)$  is asked to  be such that 
\begin{enumerate}[label=(H\arabic*)]
  \setcounter{enumi}{4}
  \item\label{VJ regularity A} $\J\in C^1(\GL_+(d))$;
  \item\label{VJ A} there exist $q>pd/(p-d)$ and $c_J>0$ such that
  \begin{equation*}
    \J(F)\geq \frac{c_J}{|\det F|^q}-\frac{1}{c_J} \quad  
    \forall  F\in \GL_+(d). \EEE
  \end{equation*}
\end{enumerate}

Finally,  for some fixed  $\delta>0$  we recall that  
$h\colon \R\rightarrow [0,1]$  is given by 
  \begin{equation}\label{eq h}
    h(\sigma)= \begin{cases}
      1 \qquad &\text{ if } \sigma\leq 0,\\
      \delta \EEE\qquad &\text{ if } \sigma>0.
    \end{cases}
\end{equation}
 In order to specify the stored elastic energy of the combined
accreting-fictitious medium, we define the functional $\W 
\colon  C(\overline{U}) \times \A\times   L^\infty( \EEE {U};  
\GL_+(d)) \rightarrow [0,\infty)
$ as 
\begin{equation*}
  \W(\sigma,y;A)\coloneqq \int_{U }h(\sigma)W(\nabla y A^{-1})+\J(\nabla y)\dx.
\end{equation*}
Here, $A\in  L^{\infty} (U;  
\GL_+(d)) $ is a placeholder for the
backstrain tensor   given by
  \eqref{eq:A0} \EEE  and $\sigma\in
C( \overline{U})$ is a placeholder for $ x \mapsto
\theta(x){-}t$,  whose sublevel set
$\{x\in U\mid \theta(x)-t<0\}$  identifies the location of the
accreting medium  at time
$t$.  
In particular,  in the accreting medium the latter energy density
reads  $W+  V^J\EEE$, whereas in the fictitious medium it is
${  \delta}W +  V^J\EEE$. Choosing $\delta$ small hence
corresponds to assuming that the fictitious material is highly elastically
compliant. Recall nonetheless that we assume that both accreting and
fictitious materials have the same $  V^J\EEE$, modeling a comparable
response to extreme compression.  

We additionally consider a second-order potential $\HH \colon W^{2,p}_{\omega}(U ;\R^d)\rightarrow [0,\infty)$ given by
\begin{equation*}
  \HH(y)\coloneqq \int_{U }H(\nabla^2 y)\dx
\end{equation*}
where $H:\R^{d\times d \times d}\rightarrow [0,\infty)$ is such that 
\begin{enumerate}[label=(H\arabic*)]
  \setcounter{enumi}{6}
  \item\label{H1 A} $H\in C^1(\R^{d\times d\times d})$ is convex;
  \item\label{H2 A} there exists a positive constant $c_H>0$ such that 
  \begin{align*}
    &c_H|G|^p-\frac{1}{c_H}\leq H(G) \leq \frac{1}{c_H}(1+|G|^p), \quad
    |D  H(G)|\leq \frac{1}{c_H}(1+|G|^{p-1}) \quad   \forall G\in \R^{d\times d\times d}\EEE\\
    &c_H|G-G'|^p\leq (DH(G)-DH(G'))\vdots (G-G') \quad   \forall
      G,\, G'\in \R^{d\times d\times d};\EEE
  \end{align*} 
  \item\label{H3 A} $H (Q G)=H (G)$ for all $G \in \R^{d\times d
      \times d}, Q \in {\rm S O}(d)$.
\end{enumerate}
Frame indifference \ref{H3 A} of $H$ is assumed to guarantee physical
consistency, albeit being not necessary for the analysis.  By
including $\mathcal{H}$ in the total energy of the combined medium, we
are indeed modeling a second-grade, nonsimple material \cite[Sec.~2.5]{Kruzik Roubicek}. The
inclusion of this second-order potential is primarily motivated by the
need to ensure sufficient compactness in the problem and the
corresponding length-scale constant $c_H$ is ideally
assumed to be very small. Moving from these considerations, the
second-order energy density $H$ is   taken to be \EEE  identical for both the
accreting and the fictitious materials, for simplicity.

\subsubsection{Viscous dissipation}
The dissipation potential  is defined via  $\RR\colon C(
\overline{U}  )\times W^{2,p}_{\omega}(U ;\R^d)\times H^{1}(U ;\R^d)\rightarrow [0,\infty)$
given by 
\begin{equation*}
  \RR(\sigma,y,\dot{y})\coloneqq \int_U  h(\sigma)R(\nabla y,\nabla \dot y)\dx 
\end{equation*}
where $R:\R^{d\times d}\times \R^{d\times d}\rightarrow [0,\infty)$ is
 specified as  
\begin{equation*}
  R(F,\dot F)\coloneqq \frac{1}{2}\dot C {:}\DD(C)\dot C \quad  
  \forall F,\, \EEE \dot F\in \R^{d\times d}
\end{equation*}
with $C\coloneqq F^\top F$ and $\dot C\coloneqq \dot F^\top F+ F^\top\dot F$.  We assume
\begin{enumerate}[label=(H\arabic*)]
  \setcounter{enumi}{9}
  \item\label{R1 A} $\DD\in C(\R^{d\times d}_{\rm sym};\R^{d\times d\times d\times d})$ is such that $\DD_{ijk\ell }=\DD_{jik\ell }=\DD_{k\ell ij}$ for every $i,j,k,\ell =1,\dots, d$; 
  \item\label{R2 A} there exists a positive constant $c_R>0$ such that 
  \begin{equation*}
    c_R|\dot C|^2\leq  \dot C {:}\DD(C)\dot C \quad   \forall C, \,
    \dot C \in \R^{d\times d}_{\rm sym}. \EEE
  \end{equation*} 
\end{enumerate}
The very structure of $R$ guarantees that it is frame indifferent \cite{Antman}.
Notice that by the definition of $R$, we have that $\partial_{\dot
  F}R$ is linear in $\dot F$.   More precisely, we have that \EEE 
\begin{equation*}
  \partial_{\dot F}R(F,\dot F)=2 F\left(\DD(C){:}\dot C\right)= 2 F\DD(F^\top F){:}(\dot{F}^\top F{+} F^\top \dot{F}).
\end{equation*}

\subsubsection{Loading and initial data}
We denote by $f\colon [0,T]\times U\rightarrow \R^d$  a given   body-force density, and we require
\begin{enumerate}[label=(H\arabic*)]
  \setcounter{enumi}{11}
  \item\label{L1 A} $f \in W^{1,\infty}\left(0,T;L^{2}(U
      ;\R^d)\right) \cap  L^\infty(Q;\R^d)$. 
\end{enumerate}

We moreover assume that the initial backstrain $A_0$ and the initial deformation $y_0$ satisfy 
\begin{enumerate}[label=(H\arabic*)]
  \setcounter{enumi}{12}
  \item\label{I1 A} $A_0\in  C(\overline{\Omega_0}; \GL_+(d))$, 
  $y_0\in \A$, and
  \begin{equation*}
    \int_U W(\nabla y_0 A_0^{-1})\mathbbm{1}_{\Omega_0}+W(\nabla y_0 )\mathbbm{1}_{U\setminus\Omega_0}+\J(\nabla y_0)+H(\nabla^2 y_0)\dx<\infty.
  \end{equation*}
\end{enumerate}
 Here and in the following,  $W(\nabla y_0
A_0^{-1})\mathbbm{1}_{\Omega_0}$ indicates the trivial extension of
$W(\nabla y_0 A_0^{-1}) $ to $U$.

\subsubsection{Growth}  Concerning the growth problem, we assume
the following 
\begin{enumerate}[label=(H\arabic*)]
  \setcounter{enumi}{13}
  \item\label{G1 A} $\gamma \in C^{0,1}(\GL_+(d)
  )$ is such that $c_\gamma\leq \gamma(\cdot)\leq C_\gamma$ for some $0<c_\gamma\leq C_\gamma$;
 \item\label{G2 A}  $\Omega_0+B_{C_\gamma T}\subset \subset U$. 
\end{enumerate}
We remark
assumption \ref{G2 A} guarantees that the accreting material   does
not \EEE reach the boundary of   the cointainer \EEE $U$ by the final time $T$, see \eqref{eq: OmegaT contained in U A} below.

\section{Notion of solution and main results}\label{sec: main results A}
\subsection{Notion of solution}
 We now specify our notion of solution to \eqref{eq growth intro
  paper 20}--\eqref{eq:cond3}.

\begin{definition}[Solution]\label{def:weak sol A}
   We say that a pair 
   \begin{equation*}
     (\theta,y)\in C^{0,1}(\overline{U}) \times  \left(L^{\infty}(0,T;W^{2,p}(U ;\R^d))\cap H^1(0,T;H^1(U ;\R^d))\right)
   \end{equation*} 
   is a \emph{solution}  to the initial-boundary-value problem \eqref{eq growth intro paper 20}--\eqref{eq:cond3}
   if
   \begin{enumerate}
   \item $\theta$ is a viscosity solution to
     \begin{align}\label{eq: eikonal system A}
       &\gamma(\nabla y(\theta(x)\wedge T,x))|\nabla \theta(x)|=1 \quad
         \text{in} \  U \setminus\overline{\Omega_0},\\
       &\theta(x)=0  \quad \text{on} \ \Omega_0, \label{eq: eikonal system A2}
     \end{align}
       namely, for any $\varphi \in C^1(\Rz^d)$, we have
     that $\gamma(\nabla y(\theta(x)\wedge T,x))|\nabla \varphi(x)|\leq
     (\geq) \,1$ at any local minimum (maximum, respectively) point
     $x\in  U \setminus\overline{\Omega_0}$ of $\varphi - \theta$. \EEE
     \item  $y(t,\cdot)\in \A$ for almost every $t\in (0,T) $;
     $y(0,\cdot)=y^0(\cdot)$ in $U$, and
     \begin{align}\label{weak sol eq A}
       &\int_0^T \!\!\!\!\int_{U }\!\left(
         h(\theta{-}t)\left( {\rm D} W(\nabla y  A^{-1})
         A^{-\top}{+}\partial_{\dot{F}} R(\nabla y, \nabla
         \dot{y})\right)\!{+} {\rm D} \J(\nabla y)\right){:}\nabla z
         \dx \dt \nonumber\\
       &\quad +\int_0^T \!\!\!\!\int_{U }\!  {\rm D} H\left(\nabla^2 y\right) {\vdots} \nabla^2 z\dx \dt =
         \int_0^T \!\!\!\!\int_{U }\! h(\theta{-}t) f\cdot z \dx \dt \notag\\[1mm]
       &\qquad \qquad \forall z\in C^{\infty}([0,T]\times \overline{U
         };\R^d) \quad \text{with}\ z = 0 \ \text{on} \ 
         [0,T]\times  ( \omega\cup \partial U)
     \end{align}
     with backstrain tensor $ A $ defined as 
     \begin{equation}\label{def A backstrain}
       A(  t, \EEE x)  \coloneqq 
       \begin{cases}
         A_0 \qquad &\text{ if } x\in  \overline{\Omega_0},\\
         \nabla y(\theta(x),x) & \text{ if } x\in  \overline{\Omega
           (  t \EEE )} \setminus \overline{\Omega_0},\\
         I & \text{ if } x\in U\setminus \overline{\Omega (  t \EEE)} ,
       \end{cases}
     \end{equation}
     where $\Omega(t):=\{x\in U \ | \ \theta(x)<t\}$ for $t\in (0,T]$.
   \end{enumerate}

\end{definition}

\subsection{Main result}

Our main result is the following.

\begin{theorem}[Existence]\label{Thm:existence sol A}
 Under assumptions {\rm\ref{W0 A}}--{\rm\ref{G2 A}}, there exists a
  solution $(\theta,y)$ to problem \eqref{eq growth intro paper 20}--\eqref{eq:cond3}.
\end{theorem}

The proof of Theorem \ref{Thm:existence sol A} is given in Section
\ref{sec:proof A}.  As already mentioned in the Introduction, Proposition
\ref{prop:mennucci A} allows us to find a solution $\theta$ to
\eqref{eq: eikonal system A}--\eqref{eq: eikonal system A2} for all
given $y\in L^\infty(0,T;W^{2,p}_\omega(U;\Rz^d))\cap H^1(0,T;H^1(U;\Rz^d))$. Then, in  Proposition \ref{prop: existence y a} we check
that, given $\theta\in C(\overline{U})$, by defining $\Omega(t)$ as
in \eqref{eq:zero}   we can find \EEE a solution $y\in L^\infty(0,T;W^{2,p}_\omega(U;\Rz^d))\cap H^1(0,T;H^1(U;\Rz^d))$ to \eqref{weak sol
  eq A}--\eqref{def A backstrain}. This allows us to implement an iterative procedure.  The
proof of  
Theorem \ref{Thm:existence sol A} follows by checking that  such iterations
converge, up to subsequence, to a   solution.   


\section{Proof of Theorem \ref{Thm:existence sol A}}\label{sec:proof A}

 We begin by recalling   \cite[Thm.~3.15]{Mennucci},   which
 ensures \EEE  the
 well-posedness of the external problem for the generalized eikonal
 equation in the whole $\Rz^d$. 
 
 \begin{proposition}[Well-posedness of the growth
    subproblem] \label{prop:mennucci A}
   Assume to
   be given $\widehat \gamma \in C(\Rz^d)$ with $c_\gamma \leq
   \widehat \gamma(\cdot) \leq C_\gamma$ for some $0<c_\gamma\leq C_\gamma$ and
   $\Omega_0\subset \Rz^d$ nonempty, open, and bounded. Then, there
   exists a unique nonnegative  continuous $\theta$ viscosity solution to
   \begin{align}
     &\widehat \gamma (x) |\nabla \theta(x)|=1 \quad \text{in} \ \ \Rz^d
     \setminus \overline{\Omega_0},\label{eq:new1 A}\\
     &\theta=0 \quad \text{on} \  \ {\Omega_0}, \label{eq:new2 A} 
   \end{align}
   namely, for any $\varphi \in C^1(\Rz^d)$, we have
     that $\widehat \gamma(x)|\nabla \varphi(x)|\leq
     (\geq) \,1$ at any local minimum (maximum, respectively) point
     $x\in  \R^d \setminus\overline{\Omega_0}$ of $\varphi - \theta$.
    Such $\theta$ is given by the representation formula 
   \begin{equation}
  \theta(x)=\min\left\{\int_0^1
    \frac{|\rho'(s)|}{ \widehat \gamma  (\rho(s))}\, \ds  \ \Big|\ \rho\in W^{1,\infty}(0,1;\R^d), \rho(0)\in \overline{\Omega_0},\rho(1)=x\right\}\!.\label{eq:representation}
\end{equation}
   In particular,  $\theta \in C^{0,1}(\Rz^d)$ with  
   \begin{equation}
   \label{eq:lip A}
   0<\frac{1}{C_\gamma} \leq |\nabla \theta(x)| \leq \frac{1}{c_\gamma} \ \
   \text{for a.e.} \ \ x\in \Rz^d.
 \end{equation}
 \end{proposition}
  The \EEE  representation formula \eqref{eq:representation}  and the bounds \ref{G1 A} ensure that
 \begin{equation}
     \frac{\operatorname{dist}(x,{\Omega_0})}{C_\gamma}\leq \theta(x)
     \leq \frac{\operatorname{dist}(x,{\Omega_0})}{c_\gamma} \quad
     \forall x \in   \Rz^d  \setminus \overline{\Omega_0}.\label{eq:boundtheta}
   \end{equation}

   Before moving to the proof of Theorem \ref{Thm:existence sol A},
   let us show that, for  all  given $\theta\in
   C(\overline{U})$, there exists a  deformation $y$   satisfying
   \eqref{weak sol eq A}--\eqref{def A backstrain}. 

\begin{proposition}[Existence for the equilibrium subproblem]\label{prop: existence y a}
 Let {\rm \ref{W0 A}}--{\rm \ref{G2
     A}} hold,  $\theta\in C(\overline{U} )$, and $\Omega(t)$ be
 defined as in \eqref{eq:zero} for all $t \in [0,T] $.  Then, there exists $y\in
 L^{\infty}(0,T;W^{2,p}_{\omega}(U ;\R^d))\cap H^1(0,T;H^1(U ;\R^d))$
 such that $y(t,\cdot)\in \A$ for almost every $t\in [0,T]$,  and
 $y(0,\cdot)=y_0(\cdot)$ in $U$, satisfying  \eqref{weak sol eq A}--\eqref{def A backstrain}. 
\end{proposition}

\begin{proof}
   We follow the blueprint \EEE of \cite{BadalFriedrichKruzik} or
 \cite{KromerRoubicek}   and argue by \EEE time-discretization.
Let $\tau\coloneqq T/{N_\tau}>0$ with $N_\tau\in \N$ given and
consider the corresponding uniform partition  of the time interval
$[0,T]$  given by   $t_i\coloneqq i \tau $, for $i=0,\dots, N_\tau$. Moreover, set $A^0_\tau\coloneqq A_0\mathbbm{1}_{\Omega_0}+I\mathbbm{1}_{U\setminus\Omega_0}$.
For $i=1,\dots, N_\tau$,   assume to know $y^j_\tau \in
\mathcal{A}$ for $j=0,1,\dots,i-1$ and define $A^i_\tau: U \to
\Rz^{d\times d}$ as  \EEE
 \begin{equation}\label{def A backstrain discrete}
    A^i_\tau(x)\coloneqq \begin{cases}
      A_0(x) \quad &\text{ if } \theta(x)=0,\\
      \nabla y^k_\tau(x) &\text{ if } \theta(x)\in(t_{k-1},t_k] \text{ for some } k=1,\dots, i{-}1,\\
      I \quad &\text{ if } \theta(x)> t_{i-1}.\\
    \end{cases}
  \end{equation}
Notice that \ref{I1 A}, the definition of $\A$, and the fact that
$p>d$, imply that  $A^i_\tau \in L^{\infty}(U;\GL_+(d))$.   We find
$y^i_\tau \in \mathcal{A}$ by solving \EEE 
  \begin{align*}
    y^i_\tau\in \argmin_{y\in \A}\left\{\W(\theta{-}t_i,y;A^i_\tau)+\HH(y)+\tau\RR\left(\theta{-}t_i,y^{i-1}_\tau,\frac{y{-}y^{i-1}_\tau}{\tau}\right) -\int_U  h(\theta{-}t_i)f\cdot y\dx \right\}.
  \end{align*} 

Under the growth conditions 
\ref{VJ A}, \ref{H2 A}, and \ref{R2 A}, the regularity and convexity assumptions \ref{W1 A}, \ref{VJ regularity A}, \ref{H1 A}, \ref{R1 A}, and \ref{L1 A}, and by using the Poincar\'e inequality \eqref{Poincare},
  the existence of $y^i_\tau\in \A$ for $i=1,\dots, N_\tau$ easily follows by the Direct Method of the calculus of variations.
Moreover, every minimizer $y_\tau^i$ satisfies the time-discrete
 Euler--Lagrange equation
\begin{align}\label{Euler-Lagrange discrete A}
  &\int_{U } h(\theta{-}t_i)\left( {\rm D}  W(\nabla y^i_\tau (A^i_\tau)^{-1})(A^i_\tau)^{- \top}{+}\partial_{\dot{F}} R\left(\nabla y^{i-1}_\tau, \frac{\nabla y^i_\tau{-}\nabla y^{i-1}_\tau}{\tau}\right)\right)\!{:}\nabla z^i  \dx \notag\\
  &+\int_U  {\rm D}  \J(\nabla y^{i}_\tau){:}\nabla z^i \dx +\int_U  {\rm D}  H\left(\nabla^2 y^i_\tau\right)\! {\vdots} \nabla^2 z^i \dx =
  \int_U  h(\theta{-}t_i)f(t_i){\cdot} z^i\dx
 \end{align}
for every $z^i\in C^{\infty}(  \overline{U} \EEE  ;\R^d)$ with $z^i\equiv 0$ on $\omega\cup \partial U$, and for every $i=1,\dots,N_\tau $.

From the minimality of $y^i_\tau$ we get that 
\begin{align*}
  \int_U &h(\theta{-}t_i)\! \left(\!W (\nabla y^i_\tau (A^i_\tau)^{-1}){+}\tau R\left(\nabla y^{i-1}_\tau,\frac{\nabla y^{i}_\tau{-\nabla y^{i-1}_\tau}}{\tau}\right)\!{-}f(t_i){\cdot} y^i_\tau\right)\!{+}\J(\nabla y^i_\tau){+}H(\nabla^2 y^i_\tau)\dx \\
  &\leq \int_U  h(\theta{-}t_{i})\left(W (\nabla y^{i-1}_\tau (A^i_\tau)^{-1})-f(t_i){\cdot} y^{i-1}_\tau\right)+\J(\nabla y^{i-1}_\tau)+H(\nabla^2 y^{i-1}_\tau)\dx\\
  &= \int_U  h(\theta{-}t_{i-1})\left(W (\nabla y^{i-1}_\tau (A^{i-1}_\tau)^{-1})-f(t_{i-1}){\cdot} y^{i-1}_\tau\right)+\J(\nabla y^{i-1}_\tau)+H(\nabla^2 y^{i-1}_\tau)\dx\\
  &\quad+\int_U (h(\theta{-}t_{i}){-}h(\theta{-}t_{i-1}))W (\nabla y^{i-1}_\tau (A^{i-1}_\tau)^{-1})\dx\\
  &\quad+\int_U h(\theta{-}t_i)\left(W(\nabla y^{i-1}_\tau (A^{i}_\tau)^{-1}){-}W(\nabla y^{i-1}_\tau (A^{i-1}_\tau)^{-1})\right) \dx\\
  &\quad -\int_U(h(\theta{-}t_i){-}h(\theta{-}t_{i-1}))f(t_i){\cdot} y^{i-1}_\tau\dx-\int_U h(\theta{-}t_{i-1})(f(t_i)-f(t_{i-1})){\cdot} y^{i-1}_\tau\dx.
\end{align*}
Summing over $i=1,\dots,n\leq N_\tau$, we obtain
\begin{align*}
  &\int_U h(\theta{-}t_n)W (\nabla y^n_\tau (A^n_\tau)^{-1}){+}\J(\nabla y^n_\tau){+}H(\nabla^2 y^n_\tau){-}h(\theta{-}t_n)f(t_n){\cdot} y^{n}_\tau\dx \\
  &\quad\quad+\sum_{i=1}^{n}\tau\int_U  h(\theta{-}t_{i})R\left(\nabla y^{i-1}_\tau,\frac{\nabla y^i_\tau{-\nabla y^{i-1}_\tau}}{\tau}\right)\!\dx\\
  &\quad\leq \int_U  h(\theta)W (\nabla y_{0} (A^{0}_\tau)^{-1})+\J(\nabla y_{0})+H(\nabla^2 y_{0})-  h(\theta)f(0){\cdot} y_{0}\dx\\
  &\quad\quad+\sum_{i=1}^{n}\int_U (h(\theta{-}t_{i}){-}h(\theta{-}t_{i-1}))W (\nabla y^{i-1}_\tau (A^{i-1}_\tau)^{-1})\dx\\
  &\quad\quad+\sum_{i=1}^{n}\int_U h(\theta{-}t_i)\left(W(\nabla y^{i-1}_\tau (A^{i}_\tau)^{-1}){-}W(\nabla y^{i-1}_\tau (A^{i-1}_\tau)^{-1})\right)\! \dx\\
  &\quad\quad -\sum_{i=1}^{n}\int_U(h(\theta{-}t_i){-}h(\theta{-}t_{i-1}))f(t_i){\cdot} y^{i-1}_\tau {+ }h(\theta{-}t_{i-1})(f(t_i)-f(t_{i-1})){\cdot} y^{i-1}_\tau\dx.
\end{align*}
The growth conditions \ref{VJ A}, \ref{H2 A}, and \ref{R2 A}, and the definition \eqref{eq h} of $h$ ensure that 
\begin{align}
  &   c_J \left\|\frac{1}{\det \nabla
     y^n_\tau}\right\|^q_{L^q(U
     )}-\frac{|U|}{c_J}+   c_H   \|\nabla^2y^n_\tau\|^p_{L^p(U ;
    \Rz^{d\times  d\times d} )} -\frac{|U|}{c_H}   \notag\\
  &\qquad   + c_{R}  \delta \EEE\sum_{i=1}^n\tau\left\|\frac{(\nabla
  y^i_\tau  - \nabla y^{i-1}_\tau)^\top}{\tau} \nabla
  y^{i-1}_\tau + (\nabla y^{i-1}_\tau)^\top \frac{\nabla
  y^i_\tau - \nabla y^{i-1}_\tau}{\tau}  \right\|^2_{L^2(U;\R^{d\times d})}   \notag\\
  &\quad \leq \int_U  h(\theta)W (\nabla y_{0} (A^{0}_\tau)^{-1})+\J(\nabla y_{0})+H(\nabla^2 y_{0})-  h(\theta)f(0){\cdot} y_{0}\dx\notag\\
  &\qquad+\sum_{i=1}^{n}\int_U (h(\theta{-}t_{i}){-}h(\theta{-}t_{i-1}))W (\nabla y^{i-1}_\tau (A^{i-1}_\tau)^{-1})\dx\notag\\
  &\qquad+\sum_{i=1}^{n}\int_U h(\theta{-}t_i)\left(W(\nabla y^{i-1}_\tau (A^{i}_\tau)^{-1}){-}W(\nabla y^{i-1}_\tau (A^{i-1}_\tau)^{-1})\right) \dx\notag\\
  &\qquad -\sum_{i=1}^{n}\int_U(h(\theta{-}t_i){-}h(\theta{-}t_{i-1}))f(t_i){\cdot} y^{i-1}_\tau {+ }h(\theta{-}t_{i-1})(f(t_i)-f(t_{i-1})){\cdot} y^{i-1}_\tau\dx.\label{discrete gronwall ineq a}
\end{align} 
 In order to obtain an a-priori estimate, we    now \EEE control the
 various terms in the above  right-hand side.  The
initial-value term is directly bounded by \ref{L1 A}--\ref{I1 A}.

 The second term in the right-hand side of \eqref{discrete
  gronwall ineq a} can be handled as follows 
\begin{align*}
  &\sum_{i=1}^{n}\int_U (h(\theta{-}t_{i}){-}h(\theta{-}t_{i-1}))W (\nabla y^{i-1}_\tau (A^{i-1}_\tau)^{-1})\dx\\
  &\quad\stackrel{\eqref{eq h}}{=}\sum_{i=1}^{n}\int_U   (1-\delta)
    \mathbbm{1}_{\{t_{i-1}< \theta\leq t_{i}\}} \EEE
    W (\nabla y^{i-1}_\tau (A^{i-1}_\tau)^{-1})\dx\\
  &\quad\leq\sum_{i=1}^{n}\int_U \mathbbm{1}_{\{t_{i-1}< \theta\leq t_{i}\}} W (\nabla y^{i-1}_\tau (A^{i-1}_\tau)^{-1})\dx\\
  &\quad=\sum_{i=1}^{n}\int_U \mathbbm{1}_{\{t_{i-1}< \theta\leq t_{i}\}} W (\nabla y^{i-1}_\tau)\dx,
\end{align*}
since $A^{i-1}_\tau(x)=I$ if $\theta(x)>t_{i-1}$.
By the growth condition \ref{W2 A}, we then have 
\begin{align}
  &\sum_{i=1}^{n}\int_U (h(\theta{-}t_{i}){-}h(\theta{-}t_{i-1}))W (\nabla y^{i-1}_\tau (A^{i-1}_\tau)^{-1})\dx\nonumber\\
  &\quad\leq \sum_{i=1}^{n}\frac{1}{c_W}(\|\nabla y^{i-1}_\tau \|^p_{L^{\infty}(U;\R^{d\times d})}+1)\int_U \mathbbm{1}_{\{t_{i-1}< \theta\leq t_{i}\}} \dx\nonumber\\
  &\quad \stackrel{\eqref{Poincare}}{\leq} c\sum_{i=1}^{n}(\|\nabla^2 y^{i-1}_\tau
    \|^p_{L^{p}(U;\R^{d\times d})}+1)\int_U \mathbbm{1}_{\{t_{i-1}<
    \theta\leq t_{i}\}} \dx\nonumber\\
  &\quad  =    c\sum_{i=1}^{n}(\|\nabla^2 y^{i-1}_\tau
    \|^p_{L^{p}(U;\R^{d\times d})}+1) \, 
    |\overline{\Omega(t_{i})}\setminus \Omega(t_{i-1})| \label{eq:coll1}
\end{align}
where  we also   the Poincar\'e inequality \eqref{Poincare} and \EEE the continuous embedding of $L^\infty(U)$
into $W^{1,p}(U)$ for $p>d$.

 As for  the third term in the right-hand side of \eqref{discrete
   gronwall ineq a}, we notice that, for $x\in \overline{\Omega_0}$,
 $A^{i}_\tau(x)=A^{i-1}_\tau(x)=A_0$. For $x\in U\setminus
 \overline{\Omega_0}$ with $\theta(x)\leq t_{i-2}$   for $i>2$,
 \EEE there exists $k\in\{1,\dots, i-2\}$ such that $\theta(x)\in
 (t_{k-1},t_k]$, and thus $A^{i}_\tau(x)=A^{i-1}_\tau(x)= \nabla
 y^k_\tau$. Similarly, for $x\in U\setminus \overline{\Omega_0}$ such
 that $\theta(x)>t_{i-1}$   for $i\geq 2 $ \EEE  we have
 $A^{i}_\tau(x)=A^{i-1}_\tau(x)=I$. Hence, the integrand is nonzero
 only for $x\in U$ such that $t_{i-2}<\theta(x)\leq t_{i-1}$    for
 $i\geq 2 $. \EEE For such $x$ we have $A^{i}_\tau(x)=\nabla y^{i-1}_\tau$ and $A^{i-1}_\tau(x)=I$, so that  
\begin{align*}
  &\sum_{i=1}^{n}\int_U h(\theta{-}t_i)\left(W(\nabla y^{i-1}_\tau (A^{i}_\tau)^{-1}){-}W(\nabla y^{i-1}_\tau (A^{i-1}_\tau)^{-1})\right) \dx\nonumber\\
  &\quad = \sum_{i=  2 \EEE}^{n}\int_U \mathbbm{1}_{\{t_{i-2}<\theta(x)\leq t_{i-1}\}}\left(W(I){-}W(\nabla y^{i-1}_\tau )\right) \dx \stackrel{\rm \ref{W2 A}}{\leq} 0.
\end{align*}

 Eventually, the regularity \ref{L1 A} ensures that
\begin{align}
  &-\sum_{i=1}^{n}\int_U(h(\theta{-}t_i){-}h(\theta{-}t_{i-1}))f(t_i){\cdot}
    y^{i-1}_\tau {+ }h(\theta{-}t_{i-1})(f(t_i)-f(t_{i-1})){\cdot}
    y^{i-1}_\tau\dx\nonumber\\
  &\quad \leq c \sum_{i=1}^{n}\int_U \mathbbm{1}_{\{t_{i-1}<
    \theta\leq t_{i}\}} |f(t_i)|\, |y^{i-1}_\tau| +
    \sum_{i=1}^{n}\tau \left\|\frac{f(t_i)-f(t_{i-1})}{\tau} \right\|_{L^2(U;\Rz^d)}
    \|y^{i-1}_\tau\|_{L^2(U;\Rz^d)}\nonumber\\
  &\quad \leq c \sum_{i=1}^{n} \| y^{i-1}_\tau\|_{L^\infty(U;\Rz^d)}
    |\overline{\Omega(t_i)}\setminus
    \Omega(t_{i-1})|\|f(t_i)\|_{L^{  \infty}(U;\Rz^d)} + c \sum_{i=1}^{n}\tau
    \|y^{i-1}_\tau\|_{L^2(U;\Rz^d)}\nonumber\\
  &\quad \stackrel{\eqref{Poincare}}{\leq}     c \sum_{i=1}^{n}
    (\|\nabla^2 y^{i-1}_\tau \|^p_{L^{p}(U;\R^{d\times d})}+1)\left(\tau +  |\overline{\Omega(t_i)}\setminus
    \Omega(t_{i-1})|\right)\label{eq:coll3}\end{align}
  where we also used the embedding $
  W^{2,p}(U;\Rz^d) \subset L^\infty(U;\Rz^d) $, the Poincar\'e inequality \eqref{Poincare}, and
  the fact that $p>2$. 

 By collecting \eqref{eq:coll1}--\eqref{eq:coll3} in
\eqref{discrete gronwall ineq a} we hence have that
\begin{align*}
   &  \|\nabla^2y^n_\tau\|^p_{L^p(U ;
    \Rz^{d\times  d\times d} )}   +\left\|\frac{1}{\det \nabla
     y^n_\tau}\right\|^q_{L^q(U
     )} \notag\\
  &\quad   + \sum_{i=1}^n\tau\left\|\frac{(\nabla
  y^i_\tau  - \nabla y^{i-1}_\tau)^\top}{\tau} \nabla
  y^{i-1}_\tau + (\nabla y^{i-1}_\tau)^\top \frac{\nabla
  y^i_\tau - \nabla y^{i-1}_\tau}{\tau}  \right\|^2_{L^2(U;\R^{d\times
    d})}   \notag\\
  &\quad \leq    c\sum_{i=1}^{n}(\|\nabla^2 y^{i-1}_\tau
    \|^p_{L^{p}(U;\R^{d\times d})}+1)\left(\tau +
    |\overline{\Omega(t_{i})}\setminus \Omega(t_{i-1})|\right) +c
\end{align*}
The Discrete  Gronwall Lemma \cite[(C.2.6), p. 534]{Kruzik
  Roubicek}  and the Poincar\'e inequality \eqref{Poincare} allow us to conclude that 
\begin{align}
  &\max_{n}\left(\|y^n_\tau\|^p_{W^{2,p}(U;\R^d)}+\left\|\frac{1}{\det \nabla
     y^n_\tau}\right\|^q_{L^q(U
     )}\right)\notag\\
  &\qquad   + \sum_{i=1}^{N_\tau}\tau \left\|\frac{(\nabla
    y^i_\tau  - \nabla y^{i-1}_\tau)^\top}{\tau} \nabla
    y^{i-1}_\tau + (\nabla y^{i-1}_\tau)^\top \frac{\nabla
    y^i_\tau - \nabla y^{i-1}_\tau}{\tau}  \right\|^2_{L^2(U;\R^{d\times d})}  \notag\\
    &\quad\leq c
    \operatorname{exp}\left(\sum_{i=1}^{N_\tau}
      |\overline{\Omega(t_{i})}\setminus \Omega(t_{i-1})|  \right)+c
    \leq c \operatorname{exp}\left(| \overline{\Omega(T)}|\right)+c\!.\label{eq:prebound A}
\end{align}

Let us now introduce the following notation for the time interpolants
of a vector $(u_{0},...,u_{N_\tau})$ over the interval $[0,T]$: We
define   the \EEE  backward-constant interpolant
$\overline{u}_{\tau}$,   the \EEE  forward-constant interpolant
$\underline{u}_{\tau}$, and   the \EEE  piecewise-affine
interpolant $ \widehat{u}_{\tau}$ \EEE on the partition $(t_i)_{i=0}^{N_\tau}$ as
  \begin{align*}
    &\overline{u}_\tau(0):=u_0, \quad \quad \overline{u}_\tau(t):=u_{i}  &\text{ if } t\in(t_{i-1},t_i] \quad\text{ for } i=1,\dots, N_\tau,\\
    &\underline{u}_{\tau}(T):=u_{N_\tau}, \quad \; \underline{u}_\tau(t):=u_{i-1} &\text{ if } t\in[t_{i-1},t_i) \quad\text{ for } i=1,\dots, N_\tau,\\
    &\widehat{u}_\tau(0):=u_0, \quad \quad \widehat{u}_\tau(t):=\frac{u_i-u_{i-1}}{t_i-t_{i-1}}(t-t_{i-1})+u_{i-1}
    &\text{ if } t\in(t_{i-1},t_i] \quad\text{ for } i=1,\dots, N_\tau.
  \end{align*}
Making use of this notation, we can rewrite \eqref{eq:prebound A} as
\begin{align}
  \!\!\|\overline{y}_\tau\|^p_{L^{\infty}(0,T;W^{2,p}(U;\R^d))}{+}\!\left\|\frac{1}{\det \nabla
     \overline{y}_\tau}\right\|^q_{L^{\infty}(0,T;L^q(U
     ))}\!\!{+}\!\!\int_0^T\!\! \| \nabla\dot{\widehat{y}}_\tau^\top
  \nabla\underline{y}_\tau {+} \nabla
  \underline{y}_\tau^\top\nabla\dot{\widehat{y}}_\tau  
  \|^2_{L^2(U;\Rz^{d{\times}d})}\!\dt  \leq c.\label{eq:prebound0 A}
\end{align}

By the Sobolev embedding of $W^{2,p}(U ;\R^d)$ into  $C^{1,1-d/p}(\overline{U}
;\R^d)$  and   by \EEE the classical result \cite[Thm. 3.1]{HealeyKromer}, the bound \eqref{eq:prebound0 A} implies
\begin{equation}\label{determinant interpolant>0 A}
   \det\nabla \overline{y}_\tau\geq     c   >0 \ \ \text{in} \ \
 [0,T]\times \overline{U}.
\end{equation}

Moreover, by the Poincar\'e inequality \eqref{Poincare}, the generalization of Korn's first
inequality by \cite{Neff} and \cite[Thm. 2.2]{Pompe}, and the  
uniform \EEE positivity of the determinant
\eqref{determinant interpolant>0 A}, it follows that  
\begin{equation*}
  \|\nabla \dot{\widehat y}_\tau\|^{   2}_{     L^2( 
    Q;\R^{d\times d} )\EEE}\leq    c \int_0^T\ \| \nabla\dot{\widehat{y}}_\tau^\top
  \nabla\underline{y}_\tau + \nabla
  \underline{y}_\tau^\top\nabla\dot{\widehat{y}}_\tau  
  \|^2_{L^2(U;\Rz^{d\times d})} \, \ds \stackrel{\eqref{eq:prebound0 A}}{\leq}c. 
\end{equation*} 
  Thus, the classical Poincar\'e inequality applied to
$\dot y$ proves that
\begin{equation}\label{bound H^1(H^1) A}
  \| {\widehat y}_\tau\|_{H^1(0,T;H^1(U; \R^d ))}\leq     c .  
\end{equation} 
Hence, the  estimates above yield 
\begin{align}
  &\overline{y}_\tau,\,\underline{y}_\tau \stackrel{*}{\rightharpoonup}
  y  \quad \text{weakly-$*$ in} \ \  L^{\infty}(0,T;W^{2,p}(U ;\R^d)),\label{convergence 1 A}\\
  &\nabla \dot{\widehat y}_\tau \rightharpoonup \nabla \dot{y}
          \quad\text{weakly in} \ \  L^{2}( Q;\R^d),\label{convergence 2 A}\\
  &\nabla \widehat y_\tau\rightarrow \nabla y \quad\text{strongly in}
    \ \  C^{0,\alpha}( \overline{Q} ;\R^d)\label{convergence 3 A} 
\end{align} 
for $\alpha\in(0, 1-d/p)$,  as
$\tau \to 0$, up to not relabeled subsequences. 
   In particular,  these convergences imply  $\det \nabla
\overline{y}_\tau \to \det \nabla y$    uniformly and, together
with the lower bound \eqref{determinant interpolant>0 A},
that $\nabla    y  
\in \GL_+(d)$ everywhere, i.e., $ y(t,\cdot) \in \mathcal{A}$ for every $t\in (0,T)$. 

Summing up the time-discrete Euler--Lagrange equations
  \eqref{Euler-Lagrange discrete A} for $i=1,\dots,N_\tau$ and rewriting in terms of the time interpolants, we get
  \begin{align}\label{Euler-Lagrange compact A} 
    &\int_0^T \int_{U } h(\theta{-}\overline{t}_\tau)\left( {\rm D
      } W(\nabla \overline{y}_\tau (\overline{A}_\tau)^{-1})(\overline{A}_\tau)^{- \top}{+}\partial_{\dot{F}} R\left(\nabla \underline{y}_\tau, \nabla \dot{\widehat{y}}_\tau\right)\right){:}\nabla \overline{z}_\tau  \dx\dt \notag
    \\
  &+\int_0^T \int_U  {\rm D}  \J(\nabla \overline{y}_\tau){:}\nabla \overline{z}_\tau +  { \rm D}  H\left(\nabla^2 \overline{y}_\tau\right) {\vdots} \nabla^2 \overline{z}_\tau \dx\dt =
  \int_0^T \int_U  h(\theta{-}\overline{t}_\tau)f(\overline{t}_\tau)\cdot \overline{z}_\tau\dx\dt.
  \end{align}

We now pass to the limit in \eqref{Euler-Lagrange compact A} as $\tau \rightarrow 0$.
Let $z\in
C^\infty(\overline{Q};\Rz^d)$ with $z\equiv 0$ on $ (0,T) \EEE\times (\omega\cup \partial U)$ be given and
let $(z^i_\tau)_{i=1}^{N_\tau}\subset W^{2,p}(U;\R^d)$ be such that $z^i_\tau\equiv 0$ on $\omega\cup \partial U$ for every $i=1,...,N_\tau$, and  $\overline{z}_\tau \to z$
strongly in $L^\infty(0,T;W^{2,p}(U;\Rz^d)) $.
First, notice that, by the coarea formula and the  Lipschitz continuity  \eqref{eq:lip A} of $\theta$, we have
\begin{equation*}
  \int_0^{\infty}\mathcal{H}^{d-1}(\partial \Omega(t))\dt=\int_U|\nabla \theta|\dx\leq \frac{|U|}{c_\gamma}<\infty.
\end{equation*}
Thus,  $|\partial \Omega(t)| =0$ for almost every $t\in
  (0,T)$. \EEE  We hence have 
\begin{equation*}
  \mathcal{L}^{d+1}\left(\left\{(t,x)\in [0,T]\times U\mid \theta(x)=t\right\}\right)=\int_0^T|\partial\Omega(t)|\dt=0,
\end{equation*}
which implies   that \EEE $h(\theta(x){-}\overline{t}_\tau(t))\to
h(\theta(x){-}t)$ for almost every $(t,x)\in Q$.
By \ref{L1 A}, it thus follows
\begin{equation*}
  \int_0^T\int_U  h(\theta{-}\overline{t}_\tau)f(\overline{t}_\tau){\cdot} \overline{z}_\tau\dx
  \dt\rightarrow \int_0^T\int_U  h(\theta{-}t)f(t){\cdot }z\dx\dt.
\end{equation*}
Similarly, for the  
dissipation, we find 
\begin{align*}
  &\!\int_0^T\!\!\!\int_U \!h(\theta{-}\overline{t}_\tau)\partial_{\dot{F}} R\left(\nabla \underline{y}_\tau,\nabla\dot{\hat y }_\tau\right){:}\nabla \overline{z}_\tau  \dx\dt\\
  &\quad=2\!\int_0^T\!\!\!\int_U \! h(\theta{-}\overline{t}_\tau)\nabla \underline{y}_\tau  \left(\DD(\nabla \underline{y}_\tau^\top \nabla \underline{y}_\tau)(\nabla\dot{\hat y }_\tau^\top \nabla \underline{y}_\tau{+}\nabla \underline{y}_\tau^\top \nabla\dot{\hat y }_\tau)\right){:}\nabla \overline{z}_\tau\dx\dt\\
  &\quad\rightarrow 2\!\int_0^T\!\!\!\int_U \! h(\theta{-}t)\nabla y  \left(\DD(\nabla y^\top \nabla y)(\nabla\dot{y}^\top \nabla y{+}\nabla y^\top \nabla\dot{y})\right){:}\nabla z\dx\dt\\
  &\quad=\!\int_0^T\!\!\!\int_U h(\theta {-}t)\partial_{\dot{F}} R\left(\nabla y,\nabla\dot{y}\right){:}\nabla z  \dx\dt
\end{align*}
by the convergences \eqref{convergence 1 A}--\eqref{convergence 3 A}, and \ref{R1 A}. 
By the continuity \ref{VJ regularity A} and the   lower bound on
the determinant \EEE \eqref{determinant interpolant>0 A},  we also have
\begin{equation*}
  \int_0^T\int_U  {\rm D}  \J (\nabla \overline{y}_\tau){:}\nabla \overline{z}_\tau\dx\dt\rightarrow \int_0^T\int_U  {\rm D}  \J (\nabla y){:}\nabla z\dx\dt. 
\end{equation*}
Moreover, convergence \eqref{convergence 3 A} guarantees that for almost every $(t,x)\in Q $,  $\overline{A}_\tau$ converges to $A$ given by \eqref{def A backstrain}.
Indeed, let $(t,x)\in Q $, $(t_{i_\tau})_\tau$   be \EEE such that
$t\in (t_{i_\tau-1},t_{i_\tau}]$ for every $\tau>0$, and
$t_{i_\tau}\rightarrow t$, as $\tau\rightarrow 0$. Thus,
$\overline{A}_\tau(t,x)=A^{i_\tau}_\tau(x)$. If $x\in
\overline{\Omega_0}$, then $A^{i_\tau}_\tau(x)=A_0(x)=A(  t, \EEE
x)$, whereas if $x\in U\setminus   \overline{\Omega(t)} \EEE$, then
$\theta(x)\geq t>t_{i_\tau-i}$ and thus, by definition \eqref{def A
  backstrain discrete}, $A^{i_\tau}_\tau(x)=I=A(  t, \EEE x)$. On
the other hand, if $x\in   \Omega(t)\EEE \setminus
\overline{\Omega_0}$, then there exists $s\in (0,t)$ such that
$\theta(x)=s$ and there exist $k_\tau\in \N$, $k_\tau\geq 1$, for
every $\tau>0$ such that $s\in (t_{k_{\tau}-1},t_{k_\tau}]$. Since
$s<t$, we can assume $t_{k_\tau}\leq t_{i_\tau-1}$, so that
$A^{i_\tau}_\tau(x)=\nabla y^{k_\tau}_\tau(x)\rightarrow \nabla y
(s,x)=\nabla y (\theta(x),x)=A(  t, \EEE x)$, by convergence \eqref{convergence 3 A}. 
Hence, by the continuity \ref{W1 A} and the bound \ref{W2 A} on $W$, convergences \eqref{convergence 1 A}--\eqref{convergence 3 A}, and dominated convergence, we have  
\begin{align*}
  \int_0^T\int_{U } h(\theta-\overline{t}_\tau)  {\rm D}  W(\nabla \overline{y}_\tau (\overline{A}_\tau)^{-1})(\overline{A}_\tau)^{- \top}{:}\nabla \overline{z}_\tau\dx\dt\\
  \rightarrow \int_0^T\int_{U } h(\theta-t)  {\rm D}  W(\theta{-}t,\nabla y A^{-1})A^{- \top}{:}\nabla z\dx\dt.
\end{align*}

The convergence of the second-gradient term follows by a standard argument 
\cite{KromerRoubicek}, which we   reproduce here \EEE for the sake of completeness. Let
$(w_{\tau}^i)_{i=1}^{N_\tau}\subset \mathcal{A}$ approximate the limiting function $y$, namely be such that  
$\overline{w}_{\tau} \to y$ strongly in
$L^\infty(0,T;W^{2,p}_\omega(U;\Rz^d))$ as $\tau\rightarrow 0$.  
Define \EEE $   \overline{z}_{\tau}
  \coloneqq
   \overline{w}_{\tau} -\overline{y}_{\tau}$.
   By convergences \eqref{convergence 1 A}--\eqref{convergence 2 A}, it follows 
that    $   \overline{z}_{\tau} \rightarrow 0$ strongly
in  $ L^\infty (0,T;H^{1}(U;\R^d))$ and $  
\overline{z}_{\tau}\stackrel{*}{ \weak }0$    weakly-$*$   in $ L^{  
  \infty } (0,T;W^{2,p} (U;\R^d))$.  
  Moreover, by the strong convergence of $\nabla^2 \overline{w}_{\tau}$ to $ \nabla^2 y$ in
  $L^{p}(   Q  ;\Rz^{d\times d \times d})$ and   by \EEE the boundedness of  ${\rm D} H (\nabla^2
  \overline{y}_\tau)$ in $L^{p'}(   Q  ;\Rz^{d\times d \times
    d})$ thanks to \ref{H3 A}, it follows 
\begin{align}
 & \limsup_{\tau\rightarrow 0}\int_0^T\!\!\int_U({ \rm D}H(\nabla^2 y)-{ \rm D} H(\nabla^2
  \overline{y}_\tau)){\vdots}(\nabla^2 y -\nabla^2 
  \overline{y}_\tau)\,\dx\,\dt\nonumber\\
  &\quad = \limsup_{\tau\rightarrow 0}\int_0^T\!\!\int_U({ \rm D}H(\nabla^2
  y)-{ \rm D} H(\nabla^2 \overline{y}_\tau)){\vdots}(\nabla^2 
    y - \nabla^2 \overline{w}_{\tau}+\nabla^2  \overline{z}_{\tau})\,\dx\,\dt\nonumber\\
  &\quad=\limsup_{\tau\rightarrow 0}\int_0^T\!\!\int_U({ \rm D}H(\nabla^2 y)-{ \rm D} H(\nabla^2 \overline{y}_\tau)){\vdots}\nabla^2  \overline{z}_{\tau}\,\dx\,\dt.\notag
\end{align}
 Hence, the Euler--Lagrange equation \eqref{Euler-Lagrange compact A} with test function $\overline{z}_\tau$   
and convergences \eqref{convergence 1 A}--\eqref{convergence 3 A}
entail   that \EEE  
\begin{align}
  &   \limsup_{\tau\rightarrow 0} \int_0^T\!\!\int_U({ \rm D}H(\nabla^2
    y)-{ \rm D} H(\nabla^2
    \overline{y}_\tau)){\vdots}(\nabla^2y - \nabla^2 \overline{y}_\tau)
     \,\dx\,\dt\nonumber\\
  &\quad=\limsup_{\tau\rightarrow 0}\Bigg(\int_0^T\!\!\int_U{ \rm D}
    H(\nabla^2 y){\vdots}\nabla^2
      \overline{z}_{\tau} +  {\rm D} \J(\nabla \overline{y}_\tau){:}\nabla \overline{z}_\tau - h(\theta{-}\overline{t}_\tau)f(\overline{t}_\tau){\cdot}   \overline{z}_{\tau} \,\dx\,\dt
  \nonumber\\
  &\qquad+\int_0^T \int_{U } h(\theta{-}\overline{t}_\tau)\left( {\rm D}  W(\nabla \overline{y}_\tau (\overline{A}_\tau)^{-1})(\overline{A}_\tau)^{- \top}{+}\partial_{\dot{F}} R\left(\nabla \underline{y}_\tau, \nabla \dot{\widehat{y}}_\tau\right)\right){:}\nabla \overline{z}_\tau \,\dx \,\dt 
    \Bigg)=0\notag
\end{align}  
   By the coercivity \ref{H2 A},  this implies that $\nabla^2
   \overline{y}_\tau \to \nabla^2y$ strongly in $L^p(Q;\Rz^{d\times d
     \times d})$. The bound on ${\rm D}H$ in \ref{H2 A}   ensures \EEE that,
   possibly passing to not relabeled subsequences,  ${\rm D}H(\nabla^2
   \overline{y}_\tau ) \weak G$ weakly in $L^{p'}(Q;\Rz^{d\times d
     \times d})$ and we can identify $G={\rm D}H(\nabla^2y)$ as $H $ is
   convex. Hence, ${ \rm D} H(\nabla^2 \overline{y}_\tau)\weak { \rm D}
  H(\nabla^2 {y})$ weakly  in $L^{p'}(Q;  
  \R^{d\times d \times d} )$ and thus \eqref{weak sol eq A} follows by
  passing to the limit in \eqref{Euler-Lagrange compact A}  taking  $\tau \to~0 $.
\end{proof}

Having checked Propositions \ref{prop:mennucci A} and \ref{prop:
  existence y a}, we proceed with the proof of Theorem
\ref{Thm:existence sol A} by   an \EEE iterative construction. We first remark that, since $y_0\in \A$ by \ref{I1 A}, 
$\nabla y_0$ is H\"older continuous and, thus, so is
the mapping $x\in U\mapsto\gamma(\nabla y_0
(   x ))$. Denoting by  $\widehat \gamma$ be any
continuous extension of such mapping to $\Rz^d$ with $c_\gamma \leq
\widehat \gamma(\cdot)\leq C_\gamma$, by
Proposition \ref{prop:mennucci A} there exists   a unique
nonnegative viscosity solution \EEE $\theta_0
\in C(\overline U )$ to problem
\begin{align*}
    &\gamma(\nabla y_0(    x ))|\nabla \theta_0(x)|=1\quad \text{in} \  \R^d  \setminus\overline{\Omega_0},\\
    &   \theta_0 =0 \quad \text{in} \  \Omega_0,
\end{align*}
satisfying \eqref{eq:lip A} in $\overline U$.
Given
$\theta=\theta_0$, on the other hand, Proposition \ref{prop: existence y a} provides the existence of $y^1
\in  L^{\infty}(0,T;W^{2,p}(U ;\R^d))\cap H^1(0,T;H^1(U
;\R^d))$ satisfying \eqref{weak sol eq A}.

For $k\geq 1$, given $y^k \in  L^{\infty}(0,T;W^{2,p}(U ;\R^d)) \cap H^1(0,T;H^1(U
  ;\R^d)) $, the map  $x\in U \mapsto \gamma(y^k(\nabla
     y^k(\theta^k(x)   {\wedge} T  ,x)))$ is H\"older continuous. 
     As above,   we \EEE extend  it continuously to $\R^d$ as $\widehat \gamma$ with $c_\gamma \leq
\widehat \gamma(\cdot)\leq C_\gamma$. Let $\theta$ be the unique nonnegative viscosity solution to \eqref{eq:new1 A}--\eqref{eq:new2 A} and let $\theta^k\in C(\overline{U})$ be its restriction to $U$. This solves 
  \begin{align*}
   &   \gamma(\nabla
     y^k(\theta^k(x)   {\wedge} T  ,x))|\nabla \theta^k (x)|=1\quad \text{in} \  U  \setminus\overline{\Omega_0},\\
      &   \theta^k   =0  \quad \text{in} \  \Omega_0
  \end{align*}
  in the   viscosity \EEE sense and fulfills \eqref{eq:lip A} in  $\overline
  U$.  Owing to the bounds \eqref{eq:boundtheta} and \ref{G2 A} we
  also have that 
   \begin{equation}\label{eq: OmegaT contained in U A}
    {\Omega}(T)=\{x\in U\mid \theta(x)<T\}\subset \subset \Omega_0+B_{C_\gamma T}\subset \subset U,  
   \end{equation}
    so that  the accreting material does not   reach \EEE the boundary of $U$ over the time interval $[0,T]$.


 For such $\theta^k$,  Proposition \ref{prop: existence y a}  applied  for
 $\theta=\theta^k$ entails the existence of a deformation $y^{k+1} \in
   L^{\infty}(0,T;W^{2,p}(U ;\R^d))\cap H^1(0,T;H^1(U ;\R^d))
  $ satisfying \eqref{weak sol eq A}.  

  The sequence  $(\theta^k,y^k )_{k\in \N}$ 
generated by this iterative process  is in general not unique, but still 
uniformly  bounded in 
$$ C^{0,1}(\overline U)\times  \left(L^{\infty}(0,T;W^{2,p}(U ;\R^d))\cap H^1(0,T;H^1(U
  ;\R^d))\right) 
 $$ 
  thanks to the bounds \eqref{eq:lip A}, \eqref{eq:prebound0 A}, and \eqref{bound H^1(H^1) A}. 
Thus, up to subsequences,
 by the Banach--Alaoglu and the Ascoli--Arzel\`a Theorems,   
  there exists   a pair  $(y,\theta)$
such that, for some $\alpha \in (0, 1-d/p)$, 
\begin{align}
  &\theta^k \to \theta \quad \text{strongly in} \  C( \overline{U}), \label{eq:conv_theta A}\\
  &y^k \stackrel{*}{\rightharpoonup} y \quad \text{weakly-$*$ in} \  L^{\infty}(0,T;W^{2,p}(U ;\R^d))\cap H^1(0,T;H^1(U  ;\R^d)),\label{eq:conv_y A}\\
  &y^k \to y \quad \text{strongly in} \  C^{1,\alpha}(
    \overline{Q}  ;\R^d), \label{eq:conv_y2 A}
\end{align}
 and $\theta$ fulfills    \eqref{eq:lip A}
 in $\overline U$.   
 By the Lipschitz continuity of $\gamma$ 
  and \EEE  by convergences 
\eqref{eq:conv_y2 A}--\eqref{eq:conv_theta A}, we  readily  have that  $x \mapsto \gamma(\nabla
  y^k(\theta^k(x)   {\wedge} T  ,x ) )$ converges to $x \mapsto \gamma(\nabla
  y(\theta(x)   {\wedge} T  ,x ) )$ uniformly in $\overline U$. Since the eikonal equation is stable with respect to the uniform convergence of
  the data \cite[Prop.~1.2]{Ishii},   $\theta$ satisfies
  \eqref{eq: eikonal system A}
  with   coefficient  $x \mapsto \gamma(\nabla
  y(\theta(x)   {\wedge} T  ,x ) )$. Moreover, since    bounds
  \eqref{eq:prebound0 A} and  \eqref{bound H^1(H^1) A} are 
  independent of $\theta$, the same arguments of the proof
of Proposition \ref{prop: existence y a} allow passing to the limit in
the Euler--Lagrange equation \eqref{weak sol eq A}, thus concluding the proof of Theorem \ref{Thm:existence sol A}.

\section{Acknowledgements}

This research was funded in
whole or in part by the Austrian Science Fund (FWF) projects 10.55776/F65,  10.55776/I5149,
10.55776/P32788,     and   10.55776/I4354,
as well as by the OeAD-WTZ project CZ 09/2023. For
open-access purposes, the authors have applied a CC BY public copyright
license to any author-accepted manuscript version arising from this
submission.

\end{document}